\documentclass[11pt]{article}

\title{A geometric proof of Jarnik's identity in the setting of weighted simultaneous approximation}
\author{Leonhard Summerer\footnote{author supported by FWF grant I 3466-N35}}

\usepackage[cp850]{inputenc}

\usepackage{amsxtra,amssymb}

\usepackage{graphicx}

\usepackage{psfrag}

\numberwithin{equation}{section}

\newtheorem{theo}{Theorem}[section]

\newtheorem{prop}[theo]{Proposition}

\newtheorem{lemm}[theo]{Lemma}

\newcommand{\bc}{{\bf {c}}}

\newcommand{\bx}{{\bf {x}}}

\newcommand{\cc}{\mathbb {C}}

\newcommand{\nn}{\mathbb {N}}

\newcommand{\qq}{\mathbb {Q}}

\newcommand{\zz}{\mathbb {Z}}

\newcommand{\rr}{\mathbb {R}}

\newcommand{\pp}{\mathbb {P}}

\newcommand{\ff}{\mathbb {F}}

\newcommand{\op}{\overline {\varphi}}

\newcommand{\up}{\underline {\varphi}}

\newcommand{\ops}{\overline {\psi}}

\newcommand{\ups}{\underline {\psi}}

\newcommand{\of}{\overline {\phi}}

\newcommand{\uf}{\underline {\phi}}

\newcommand{\ti}{\underline {t}}

\newcommand{\Ti}{\overline {t}}

\newcommand{\si}{\underline {s}}

\newcommand{\bw}{{\bf {w}}}

\begin{document}

\maketitle

\begin{abstract} Jarnik's identity plays a major role in classical simultaneous approximation to two real numbers. Recently O. German [2] has shown a generalization to the weighted setting in which the identity has to be replaced by two inequalities. His methods belong to classical geometry of numbers. The aim of this paper is to provide an alternative approach based on a careful examination of certain successive minima functions that stem from parametric geometry of numbers, a method that has already been successfully employed to generalize Jarnik's identity to higher dimensions in the classical setup in [3] and [7]. 
\vspace{6mm}

 
\end{abstract}

\section{ Weighted simultaneous approximation }

Simultaneous approximation to $m$ real numbers $\xi_1,\ldots,\xi_m$ with $1,\xi_1,\ldots,\xi_{m}$ linearly independent over $\qq$ with respect to a system of exponents $(1,-\nu_1,\ldots,-\nu_m)$ with 
$$\nu_i>0 \mbox{ for } 1\leq i\leq m\mbox{ and }\nu_1+\ldots+\nu_m=1$$
leads to consider non-trivial solutions $\bx:=(x,y_1,\ldots,y_{m})\in\zz^{m+1}$ of the system
\begin{eqnarray*}
|x|&\leq& e^q\\
|\xi_1x-y_1|&\leq&e^{-\nu_1q}\\
\vdots&&\vdots\\
|\xi_{m}x-y_{m}|&\leq&e^{-\nu_{m}q} 
\end{eqnarray*} 
for any parameter $q>0$. If ${\mathcal B}(q)$ consists of points $(p_0,p_1,\ldots,p_{m})$ with $|p_0|\leq e^q$, $|p_i|\leq e^{-\nu_i q}$ for $1\leq i\leq m$, and $\Lambda=\Lambda({\bf \xi})$ the lattice of points $p(\bx):=(x,\xi_1x-y_1,\ldots,\xi_{m}x-y_{m})$ with $(x,y_1,\ldots,y_{m})\in\zz^{m+1}$, Minkowski's Theorem guarantees that there is a nonzero lattice point in ${\mathcal B}(q)$, in other words the first minimum $\lambda_1(q)$ with respect to ${\mathcal B}(q)$ and $\Lambda$ is at most $1$. It has turned out that in the case of classical simultaneous approximation, i.e. with exponents $(1,-1/m,\ldots,-1/m)$, the study of the successive minima functions $\lambda_1(q),\ldots,\lambda_{m+1}(q)$ provides a useful tool for establishing relations between the approximation constants attached to $\xi_1,\ldots,\xi_m$. In particular the famous Jarnik identity in the case $m=2$ of approximation to two reals (see [6]) and its version for higher dimensions as worked out by W.M. Schmidt and the author in [7] as well as independently by A. Marnat in [3], can be proved by a closer examination of the joint behaviour of $\lambda_1(q),\lambda_2(q),\lambda_3(q)$. 
\vspace{4mm}

The goal of the following investigations is a generalization of this approach in dimension two to the case of a system of exponents $(1,-\nu_1,-\nu_2)$ where
$$\nu_1\leq\nu_2 \mbox{ and } \nu_1+\nu_2=1.\leqno(1.1)$$
They were motivated by [5] and the generalization of Khintchine's transference principle to the weighted setting in [1]. O. German has obtained a very general transference principle for the weighted setting in terms of weighted diophantine exponents that generalize the classical diophantine exponents $\omega,\hat\omega,\omega^*,\hat\omega^*$. Specialising to the case of two reals, he states the weighted version of Jarnik's identity as Theorem 5, page 4 in [2]. In this paper we will introduce approximation constants inspired by the geometric approach and deduce a similar result. The more delicate problem of determining the joint spectrum of the approximation constants in this case, as determined by Roy [4] in the classical case is still open.  
\vspace{4mm}

As in [6] for given $\bx\in\zz^3\setminus\{{\bf{0}}\}$ we let $\lambda_{\bx}(q)$ be the least $\lambda>0$ with $p(\bx)\in\lambda{\mathcal B}(q)$. It follows that
$$\lambda_{\bx}(q)=\max\{|x|e^{-q},|\xi_1x-y_1|e^{\nu_1 q},|\xi_2x-y_2|e^{\nu_2 q}\}$$
and for $1\leq i\leq 3$ we have $\lambda_i(q)=\lambda_{\bx}(q)$ for some $\bx$ depending on $q$ and $i$. Rather than with these functions we will work with their logarithms 
$$L_i(q):=\log\lambda_i(q)\mbox{ and } L_{\bx}:=\log \lambda_{\bx}(q)$$
since this definition implies that $L_{\bx}(q)$ is piecewise linear with slopes $-1,\nu_1,\nu_2$ with the additional condition that for fixed ${\bx}$ the slope of $L_{\bx}$ can only increase with $q$. As a consequence, the functions $L_i(q)$ are continuous and piecewise linear with slopes from $\{-1,\nu_1,\nu_2\}$. Hence for $i=1,2,3$ the quantities
$$\up_i^{\bw}:=\liminf \frac{L_i(q)}{q}\mbox{ and }\op_i^{\bw}:=\limsup \frac{L_i(q)}{q}\leqno(1.2)$$
are finite; they are intended to quantify the joint approximability of $(\xi_1,\xi_2)$ with weight $\bw:=(\nu_1,\nu_2)$. Moreover by Minkowski's Theorem we get
$$|L_1(q)+L_2(q)+L_3(q)|\leq c\leqno(1.3)$$
with some absolute constant $c>0$. 

In [6] it is explained in detail how the triple of functions $(L_1,L_2,L_3)$ can be replaced in a canonical way by another triple $(P_1,P_2,P_3)$ that has the property 
$$P_1(q)+P_2(q)+P_3(q)=0\leqno(1.4)$$
for any $q>0$ and satisfies $|P_i(q)-L_i(q)|<2c$ for $i=1,2,3$. The slopes of $P_1$ and $P_3$ are still among $-1,\nu_1,\nu_2$ but $P_2$ may in addition have slopes $2,-2\nu_1,-2\nu_2$ on intervals of length less than $4c$. It is easily deduced that in an interval where $P_3(q)-P_2(q)>4c$ the function $P_3$ has no local minimum and that in an interval where $P_2(q)-P_1(q)>4c$ the function $P_1$ has no local maximum. If $\{i,j,k\}=\{1,2,3\}$ and $P_i$ has slope $-1$ in some interval, then 
$$|P_j(q)-P_j(q')-\nu_1(q-q')|<4c\mbox{ and }|P_k(q)-P_k(q')-\nu_2(q-q')|<4c$$
or vice versa for $q,q'$ in that interval.

\section{The global behaviour of $(P_1,P_2,P_3)$}

We introduce the functions $\psi_i(q):=\frac{P_i(q)}{q}$ for $i=1,2,3$. Then 
$$\psi_i(q)-\frac{L_i(q)}{q}\ll1/q$$
so that the definition in $(1.2)$ yields
$$\op_i^{\bw}=\limsup \psi_i(q)\mbox{ and }\up_i^{\bw}=\liminf \psi_i(q).\leqno{(3.1)}$$
Moreover, quite similar to section 7, p. 86 of [6], there exist functions $g_1,g_3$ tending to infinity such that
$$\psi_1(q)\geq-1+\frac{g_1(e^q)}{q}\mbox{ and }\psi_3(q)\leq\nu_2-\frac{g_3(e^q)}{q}\leqno{(3.2)}$$
so that
$$-1<\psi_i(q)<\nu_2\mbox{ for } i=1,2,3.$$

If we have 
$$P_3(q)-P_1(q)\leq C\leqno(2.1)$$
for some constant $C$ to be specified below and all large $q$, then $\psi_i(q)$ tends to zero and hence $\up_i^{\bw}=\op_i^{\bw}=0$ for $i=1,2,3$. If all large $q$ lie in a sequence of intervals which alternate between $q$'s for which condition $(2.1)$ is satisfied and $q$'s for which it isn't, then still $\op_1^{\bw}=\up_3^{\bw}=0$. Hence we will be interested in the case where   
$$P_3(q)-P_1(q)>C\leqno(2.2)$$
for all large $q$.

By the linear independence of $1,\xi_1,\xi_{2}$ there are arbitrarily large values of $p$ with $L_3(p)=L_2(p)$ (see Corollary 2.2 in [6], the paper [8] is devoted to an alternative criterion in the case of linearly dependent reals ), hence $P_3(p)-P_2(p)<4c$, and arbitrarily large values of $p^*$ with $L_2(p^*)=L_1(p^*)$, hence $P_2(p^*)-P_1(p^*)<4c$. If $C$ is chosen greater than $8c$, $P_3(p)-P_2(p)<4c$ implies $P_2(p)-P_1(p)>4c$ and similarly $P_2(p^*)-P_1(p^*)<4c$ implies $P_3(p^*)-P_2(p^*)>4c$. So by the intermediate value Theorem there are arbitrarily large values of $p$ with
$$P_3(p)-P_2(p)=4c,$$
as well as arbitrarily large numbers $p^*$ with
$$P_2(p^*)-P_1(p^*)=4c.$$
When $p,p^*$ are such numbers, then $P_3(p^*)-P_2(p^*)>C-4c$ and 
$$P_3(p^*)-P_2(p^*)-(P_3(p)-P_2(p))\leq 3|p-p^*|,\leqno(2.3)$$
so that $P_3(p)-P_2(p)>C-4c-3|p-p^*|$, which together with $P_3(p)-P_2(p)=4c$ yields $3|p-p^*|>C-8c$, hence $|p-p^*|>C/3-8c/3>4c$ if $C>20c$.

For every $p$ with $P_3(p)-P_2(p)=4c$, hence $P_2(p)-P_1(p)>4c$, there is a smallest $a$ and a largest $b$ with $a\leq p\leq b$ such that
$$P_3(a)-P_2(a)=P_3(b)-P_2(b)=4c$$
and $P_2(q)-P_1(q)>\gamma$ for $a\leq q\leq b$. Such an interval $[a,b]=:T$ will be called a \emph{top interval}. It is not required that $P_3(q)-P_2(q)\leq4c$ in this interval. Also it may happen that $a=p=b$, so that the interval consists of a single number.\\
For every $p^*$ with $P_2(p^*)-P_1(p^*)=4c$ there is a smallest $a^*$ and a largest $b^*$ with $a^*\leq p^*\leq b^*$ such that
$$P_2(a^*)-P_1(a^*)=P_2(b^*)-P_1(b^*)=4c$$
and $P_3(q)-P_2(q)>4c$ for $a^*\leq q\leq b^*$. Such an interval $[a^*,b^*]=:B$ will be called a \emph{bottom interval}. By $(2.3)$, a top interval has distance greater than $4c$ from a bottom interval.

For intervals $[r,s]$ and $[r',s']$ we write $[r,s]<[r',s']$ if $s<r'$. We may thus arrange all the top resp. bottom intervals into a sequence 
$$T_1<T_2<T_3<\ldots,\mbox{ resp. }B_1<B_2<B_3<\ldots.$$
There cannot be two adjacent top intervals in this sequence: for if $[a,b]<[a',b']$ were two such intervals, then $ P_2(q)-P_1(q)>4c$ for $a\leq q\leq b'$, and $b$ would no longer be the largest number as required in the definition of top intervals. Similarly, there cannot be two adjacent bottom intervals. Hence our sequence becomes
$$T_1<B_1<T_2<B_2<\ldots,$$
where each $T_j$ is a top interval and each $B_j$ a bottom interval. If $T_j=[a_j,b_j]$ and $B_j=[a^*_j,b^*_j]$ we have
$$\ldots<a^*_{j-1}\leq b^*_{j-1}<a_j\leq b_j<a^*_j\leq b^*_j<a_{j+1}\leq b_{j+1}<\ldots.$$

For $q$ in $(b^*_{j-1},a^*_j)$ we have $P_2(q)-P_1(q)>4c$ and hence $L_2(q)\neq L_1(q)$ which implies that the function $P_1$ has no local maximum (see [6], section 5, p. 81). There will be some $p^*_j$ in $[b^*_{j-1},a^*_j]$ such that $P_1$ is decreasing for $b^*_{j-1}\leq q\leq p^*_j$ and increasing for $p^*_j\leq q\leq a^*_j$ or $P_1$ is increasing in the whole interval. $P_1$ cannot be decreasing in the whole interval for if this were so, $P_2-P_1$ would be increasing on this interval, a contradiction to $P_2(a^*_j)-P_1(a^*_j)=4c$. 

 Also there will be a $p_j$ in $[b_j,a_{j+1}]$ such that $P_3$ is increasing for $b_j\leq q\leq p_j$ and decreasing for $p_j\leq q\leq a_{j+1}$ or $P_3$ is increasing in the whole interval. $P_3$ cannot be decreasing in the whole interval as $P_3(b_j)-P_2(b_j)=4c$ and the difference $P_3-P_2$ would not increase if $P_3$ was decreasing to the right of $b_j$.  

When $P_1$ has a minimum in $p^*_{j+1}\in[b^*_{j},a^*_{j+1}]$ we call $B_{j}$ a bottom interval of type 1. If $P_1$ is increasing in $[b^*_{j},a^*_{j+1}]$ we call $B_{j}$ a bottom interval of type 2. Analoguously, if $P_3$ has a maximum in $p_{j-1}\in[b_{j-1},a_{j}]$ we call $T_{j}$ a top interval of type 1 and if $P_3$ is increasing in $[b_{j-1},a_{j}]$ we call $T_{j}$ a top interval of type 2.

Let $j$ be fixed now and assume $j$ is large, so that $q\in[b^*_{j-1},a_{j+1}]$ is large. The interval $[b^*_{j-1},a_{j+1}]$ starts at the end of $B_{j-1}$ and ends at the beginning of $T_{j+1}$, hence contains $T_j$ and $B_j$ as well as the interval $I_j:=[b_j,a_j^*]$ lying between $T_j$ and $B_j$. 
In $I_j$ we have $P_2(q)-P_1(q)>4c$ and $P_3(q)-P_2(q)>4c$, hence $P_1$ has no local maximum and $P_3$ has no local minimum within $I_j$. We claim that both functions must be increasing in this interval. If this were not so, we would have $p^*_j>b_j$ so that 
$$P_1\mbox{ is decreasing for } b_{j}\leq q\leq p^*_j\mbox{ and increasing for }p^*_j\leq q\leq a^*_j,\leqno(2.4)$$ 
or $p_j<a^*_j$ such that 
$$P_3\mbox{ is increasing for } b_j\leq q\leq p_j\mbox{ and decreasing for }p_j\leq q\leq a^*_{j}.\leqno(2.5)$$
If $(2.4)$ holds and $p=\min\{p_j,p^*_j\}$ it follows that $P_3$ is increasing in $[b_j,p]$ with slope $\nu_1$ resp. $\nu_2$ and $P_1$ decreasing in $[b_j,p^*_j]$ with slope $-1$, so that $P_2$ is increasing with slope $\nu_2$ resp. $\nu_1$. Hence 
$$P_3(q)-P_2(q)=P_3(b_j)-P_2(b_j)=4c$$
in $[b_j,p]$, and there is no $q$ in this interval with $P_2(q)-P_1(q)\leq4c$. This contradicts the maximality property of the right endpoint $b_j$ of $T_j=[a_j,b_j]$ if $p^*_j>b_j$.\\
If $(2.5)$ holds, in an analoguous manner we obtain a contradiction to the minimality property of the left endpoint $a^*_j$ of $B_j=[a^*_j,b^*_j]$.

Consequently $P_2$ decreases with slope $-1$ in $I_j$ and regarding the slopes of $P_1$ and $P_3$ we claim that there exists some $r_j\in I_j$ such that $P_1(q)$ increases with slope $\nu_1$ in $[b_j,r_j]$ and with slope $\nu_2$ in $[r_j,a^*_j]$. In turn, $P_1(q)$ increases with slope $\nu_2$ in $[b_j,r_j]$ and with slope $\nu_1$ in $[r_j,a^*_j]$ (see Figure 1 below for the position of $r_j$). Note that the cases $r_j=b_j$ and $r_j=a^*_j$ are not ruled out in which case the slopes of $P_1$ resp. $P_3$ are constant in $I_j$.

In fact, if $P_1$ changes slope between $\nu_1$ and $\nu_2$ at some $r_j\in I_j$, the fact that $P_2$ doesn't together with $(1.3)$ implies that $P_3$ changes its slope in the opposite way. However, for $P_1$ such a change of slope can only increase the slope of $P_1$ as $P_2(q)-P_1(q)>4c$, hence $L_2(q)\neq L_1(q)$ for $q\in I_j$ so that $L_1(q)=L_{\bx}(q)$ for the same $\bx$ for all $q\in I_j$. 

The following picture illustrates the possible behaviour of $P_1,P_2,P_3$ for $p^*_j\leq q\leq b_{j+1}$ in the case of two successive top intervals $T_j$ and $T_{j+1}$ of type 1 with a bottom interval $B_j$ of type 1 between them. Within top resp. bottom intervals $(P_3+P_2)/2$ resp. $(P_1+P_2)/2$ are indicated by dotted lines.

\vspace{5mm}
\scalebox{0.5}{
\psfrag{p}{\Large$p_{j-1}=p^*_j$}
\psfrag{s}{\Large$p_j=p^*_{j+1}$}
\psfrag{r}{\Large$r_j$}
\psfrag{a}{\Large$a_j$}
\psfrag{b}{\Large$b_j$}
\psfrag{c}{\Large$a^*_j$}
\psfrag{d}{\Large$b^*_j$}
\psfrag{e}{\Large$a_{j+1}$}
\psfrag{f}{\Large$b_{j+1}$}
\psfrag{F}{\Huge Figure 1}

\includegraphics{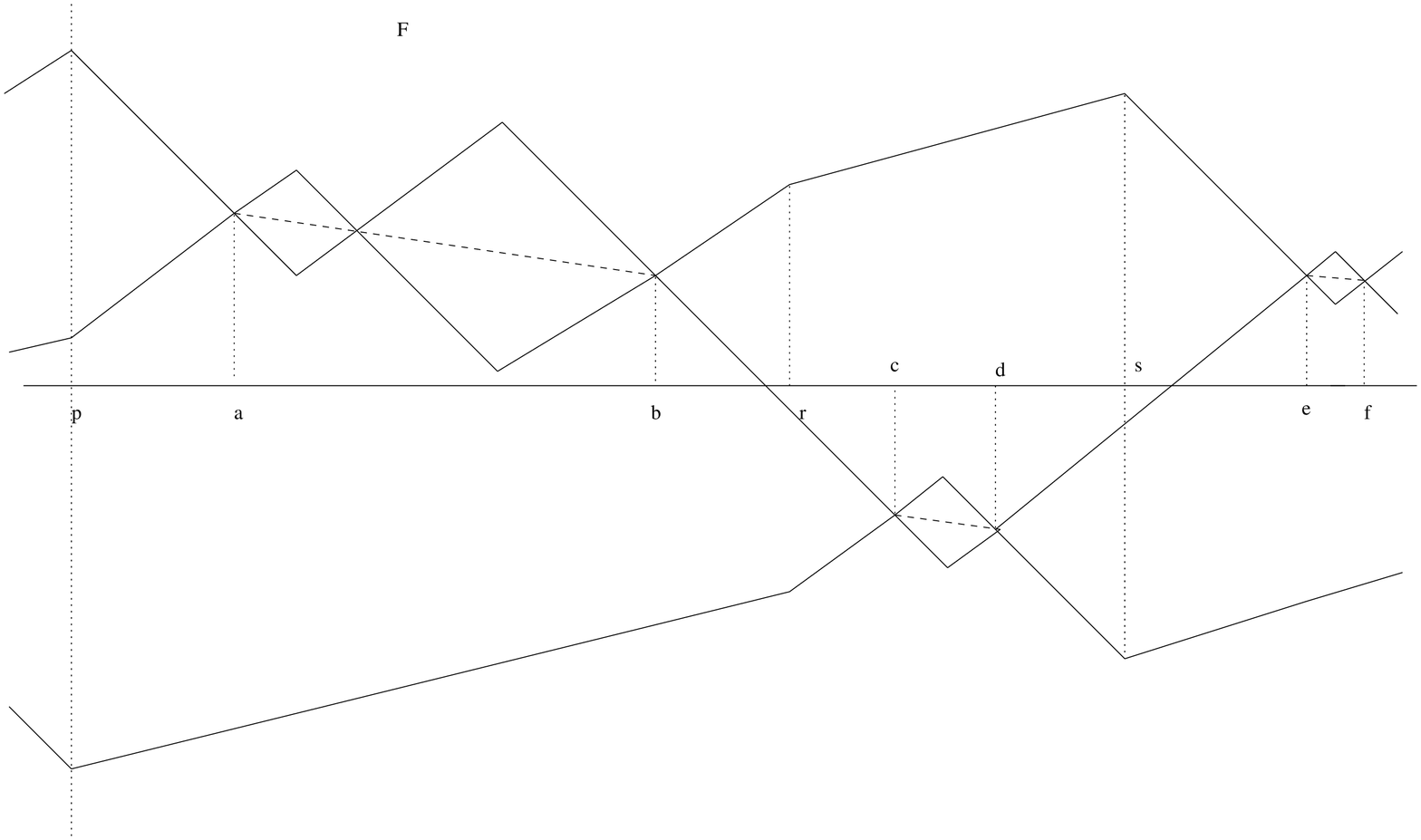}}
\vspace{4mm}

A possible behaviour of $P_1,P_2,P_3$ in case of a top interval $T_j$ of type 1 followed by a bottom interval $B_j$ of type 2 and a top Interval $T_{j+1}$ of type 1 looks as follows: 
\vspace{4mm} 

\scalebox{0.5}{
\psfrag{p}{\Large$p_{j-1}=p^*_j$}
\psfrag{q}{\Large$p_j$}
\psfrag{r}{\Large$r_j$}
\psfrag{a}{\Large$a_j$}
\psfrag{b}{\Large$b_j$}
\psfrag{c}{\Large$a^*_j$}
\psfrag{d}{\Large$b^*_j$}
\psfrag{e}{\Large$\!\!a_{j+1}$}
\psfrag{f}{\Large$\!b_{j+1}$}
\psfrag{F}{\Huge Figure 2}

\includegraphics{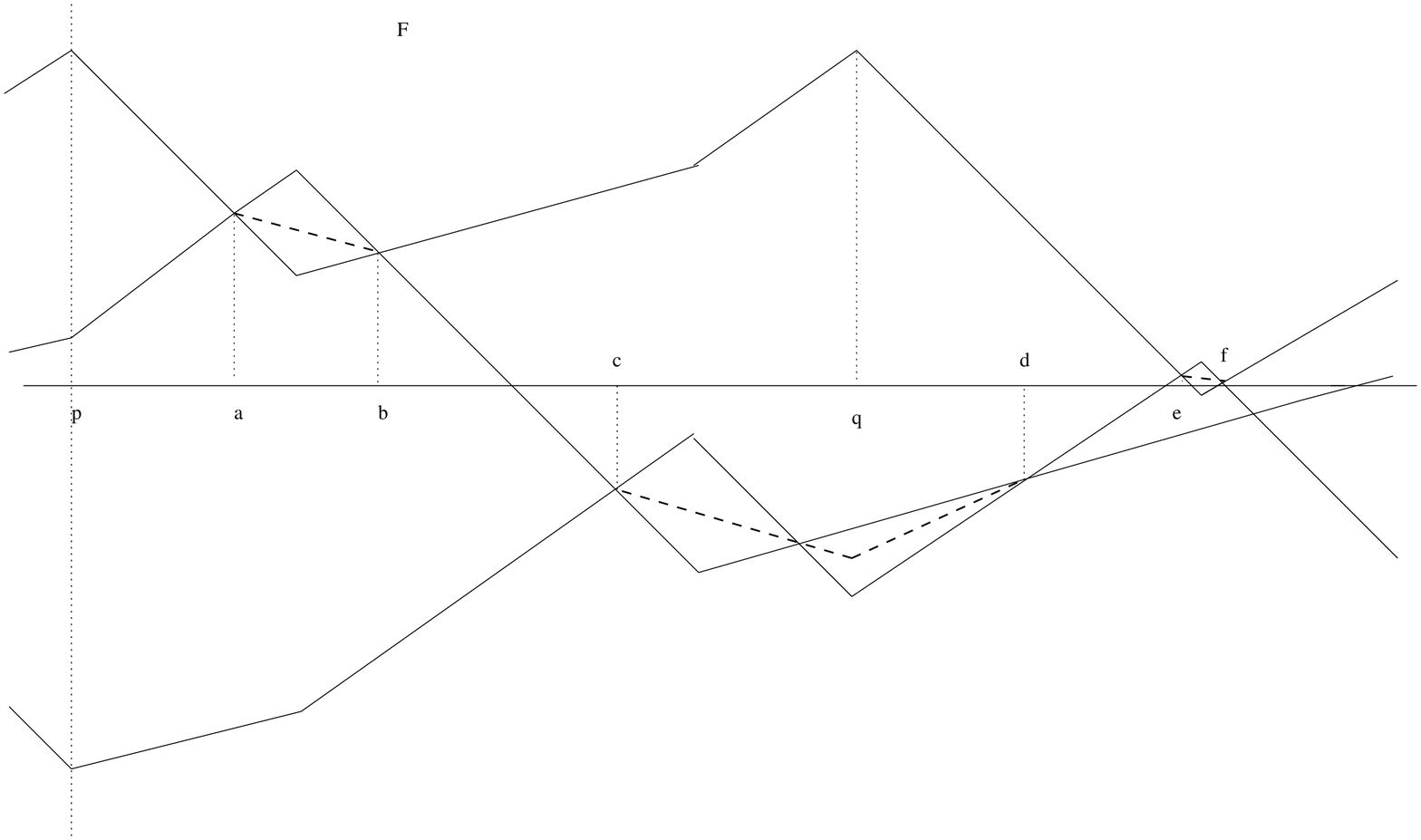}}
\vspace{4mm}

Finally, an example of a graph showing a top interval $T_j$ of type 1 followed by a bottom interval $B_j$ of type 1 and a top Interval $T_{j+1}$ of type 2:

\vspace{4mm}
\scalebox{0.5}{
\psfrag{F}{\Huge Figure 3}
\psfrag{p}{\Large$p^*_j$}
\psfrag{q}{\Large$p_{j-1}$}
\psfrag{r}{\Large$p^*_{j+1}$}
\psfrag{f}{\Large$b_{j+1}$}
\psfrag{a}{\Large$a_j$}
\psfrag{b}{\Large$b_j$}
\psfrag{c}{\Large$a^*_j$}
\psfrag{d}{\Large$b^*_j$}
\psfrag{e}{\Large$a_{j+1}$}

\includegraphics{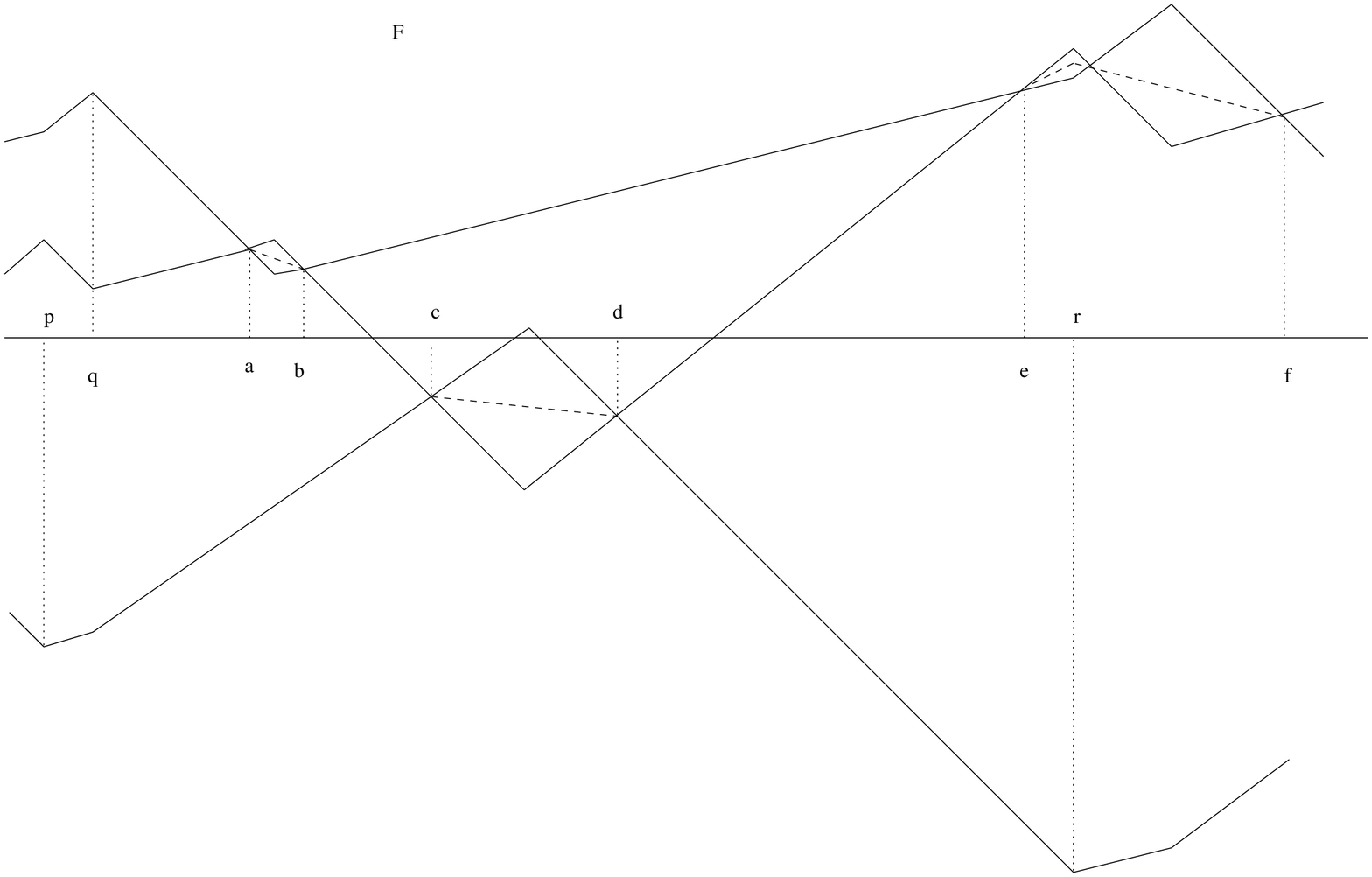}}
\vspace{4mm}

\section{The behaviour of $\psi_1(q)$ and $\psi_3(q)$}

\begin{lemm}\label{L+L} Cosider a fixed, large $j$ and drop the subscript $j$ in $I_j$ and all the points lying in this interval. Then 
$$(2-\nu_1)\psi_1(a^*)+(1+\nu_1)\psi_3(b)-O(1/b)\leq-3\psi_1(a^*)\psi_3(b)\leq(2-\nu_2)\psi_1(a^*)+(1+\nu_2)\psi_3(b)+O(1/b),$$

where the implied constants depend only on $c$, hence are absolute.
\end{lemm}

Proof: We have $P_2(a^*)-P_1(a^*)=O(1)$, hence $2P_1(a^*)+P_3(a^*)=O(1)$. Similarly, $2P_3(b)+P_1(b)=O(1)$. $P_3$ is increasing in $[b,r]$ with slope $\nu_2$ and with slope $\nu_1$ in $[r,a^*]$, so that
$$P_3(a^*)-P_3(b)=\nu_1(a^*-r)+\nu_2(r-b)=:\nu(a^*-b)$$
for some $\nu$ with $\nu_1\leq\nu\leq\nu_2$ and hence
$$2P_1(a^*)+P_3(b)=-\nu(a^*-b)+O(1).\leqno(3.3)$$
$P_2$ is decreasing in $I$ with slope $-1$ and $P_3(b)-P_2(b)=O(1)$ as well as $P_2(a^*)-P_1(a^*)=O(1)$ so that
$$P_3(b)-P_1(a^*)=a^*-b+O(1).\leqno(3.4)$$
Multiplying $(3.4)$ by two and adding it to $(3.3)$ yields
$$3P_3(b)=(2-\nu)(a^*-b)+O(1),$$
and substituting $\frac{3}{2-\nu}P_3(b)$ for $a^*-b$ in $(3.4)$ gives
$$P_1(a^*)=-\frac{1+\nu}{2-\nu}P_3(b)+O(1).\leqno(3.5)$$ 

Therefore
\begin{eqnarray*}
\psi_1(a^*)+\frac{1+\nu}{2-\nu}\psi_3(b)&=&\frac{1}{a^*b}\left(bP_1(a^*)+a^*\frac{1+\nu}{2-\nu}P_3(b)\right)\\
&=&\frac{1}{a^*b}\left(bP_1(a^*)-a^*P_1(a^*)\right)+O(1/b)\\
&=&\frac{1}{a^*b}\left((b-a^*)P_1(a^*)\right)+O(1/b)\\
&=&-\frac{3P_3(b)}{(2-\nu)a^*b}P_1(a^*)+O(1/b)\\
&=&-\frac{3}{2-\nu}\psi_1(a^*)\psi_3(b)+O(1/b),
\end{eqnarray*}
which we can rewrite as
$$(2-\nu)\psi_1(a^*)+(1+\nu)\psi_3(b)=-3\psi_1(a^*)\psi_3(b)+O(1/b).\leqno(3.6)$$
With $\nu_1\leq\nu\leq\nu_2$ and taking into account that $\psi_1(a^*)\leq0$ and $\psi_3(b)\geq0$ our claim easily follows.
\vspace{4mm}

\begin{lemm}\label{LI}
(i) If $P_1$ is decreasing in an interval with large end points, then also $\psi_1$ is decreasing in that interval. If $P_1$ is increasing in an interval with large end points, then also $\psi_1$ is increasing in that interval.

(ii) For $q\in B_j=[a^*_j,b^*_j]$ with large $j$ we have $\psi_1(q)\leq\max\{\psi_1(a^*_j),\psi_1(a^*_{j+1})\}+O(1/a^*_j)$ provided $\psi_1(a^*_j)>-\nu_1/2$ and $\psi_1(a^*_ {j+1})>-\nu_1$.
\end{lemm}

Proof: (i) Suppose $q>q'$ are in the interval. When $P_1$ is decreasing, then  $P_1$ has slope $-1$, hence $P_1(q)-P_1(q')=-(q-q')$, so that
$$\psi_1(q)=\psi_1(q')q'/q-(1-q'/q)=\psi_1(q')-(1+\psi_1(q'))(1-q'/q)<\psi_1(q')$$
by $(3.2)$.\\
When $P_1$ is increasing, then $P_1(q)-P_1(q')\geq\nu_1(q-q')$, yielding
$$\psi_1(q)\geq\psi_1(q')q'/q+\nu_1(1-q'/q)=\psi_1(q')+(\nu_1-\psi_1(q'))(1-q'/q)>\psi_1(q')$$
as $\psi_1(q')\leq0$ in view of $(1.4)$. 

(ii) If $B_j$ is of type 1 we have
\begin{eqnarray*}
P_1(q)&\leq&\frac{1}{2}\left(P_1(q)+P_2(q)\right)+O(1)\\
&\leq&\frac{1}{2}\left(P_1(a^*_j)+P_2(a^*_j)\right)+\frac{\nu_2-1}{2}(q-a^*_j)+O(1)\\
&=&P_1(a^*_j)+\frac{\nu_2-1}{2}(q-a^*_j)+O(1)\\
&=&P_1(a^*_j)-\frac{\nu_1}{2}(q-a^*_j)+O(1),
\end{eqnarray*}
hence the assumption $\psi_1(a^*_j)>-\frac{\nu_1}{2}$ yields
\begin{eqnarray*}
\psi_1(q)&\leq&\psi_1(a^*_j)a^*_j/q-\frac{\nu_1}{2}(1-a^*_j/q)+O(1/q)\\
&=&\psi_1(a^*_j)-(\frac{\nu_1}{2}+\psi_1(a^*_j))(1-a^*_j/q)+O(1/q)\\
&\leq&\psi_1(a^*_j)+O(1/a^*_j).
\end{eqnarray*}

If $B_j$ is of type 2 then 
\begin{eqnarray*}
P_1(q)&\leq&\frac{1}{2}\left(P_1(q)+P_2(q)\right)+O(1)\\
&\leq&\max\left\{\!\frac{1}{2}\!\left(P_1(a^*_j)\!+\!P_2(a^*_j)\right)\!-\!\frac{1\!-\!\nu_2}{2}(q-a^*_j),\frac{1}{2}\left(P_1(b^*_j)\!+\!P_2(b^*_j)\right)\!-\!\frac{1}{2}\!(b^*_j-q)\!\right\}\!+\!O(1)\\
&=&\max\left\{P_1(a^*_j)-\frac{1-\nu_2}{2}(q-a^*_j),P_1(b^*_j)-\frac{1}{2}(b^*_j-q)\right\}+O(1)\\
&=&\max\left\{P_1(a^*_j)-\frac{\nu_1}{2}(q-a^*_j),P_1(b^*_j)-\frac{1}{2}(b^*_j-q)\right\}+O(1).
\end{eqnarray*}
If this maximum is $P_1(a^*_j)-\frac{\nu_1}{2}(q-a^*_j)$, then we conclude that $\psi_1(q)\leq\psi_1(a^*_{j})+O(1/a^*_{j})$ as for type 1 intervals. If the minimum is $P_1(b^*_j)-\frac{1}{2}(b^*_j-q)$, the fact that $P_1$ is increasing in $[b^*_{j},a^*_{j+1}]$ yields $P_1(b^*_j)\leq P_1(a^*_{j+1})-\nu_1(a^*_ {j+1}-b^*_{j})$ and hence $\nu_1<1/2$ implies
$$P_1(q)\leq P_1(a^*_{j+1})-\nu_1(a^*_{j+1}-q)+O(1).\leqno(3.8)$$
If $\psi_1(a^*_j)>-\nu_1$ we may thus write
\begin{eqnarray*}
\psi_1(q)&\leq&\psi_1(a^*_{j+1})a^*_{j+1}/q-\nu_1(a^*_{j+1}/q-1)+O(1/q)\\
&=&\psi_1(a^*_{j+1})-(\nu_1+\psi_1(a^*_{j+1}))(a^*_{j+1}/q-1)+O(1/q)\\
&\leq&\psi_1(a^*_{j+1})+O(1/a^*_{j}).
\end{eqnarray*}

\begin{lemm}\label{LII}
(i) If $P_3$ is decreasing in an interval with large end points, then also $\psi_3$ is decreasing in that interval. If $P_3$ is increasing with slope $\nu_2$ in an interval with large end points then also $\psi_3$ is increasing in that interval. If $P_3$ is increasing with slope $\nu_1$ in an interval with large end points then $\psi_3$ is increasing in that interval provided 
$$\psi_3(q)<\nu_1\mbox{ for all }q \mbox{ in that interval}.\leqno(3.7)$$

(ii) For $q\in T_j=[a_j,b_j]$ with large $j$ we have $\psi_3(q)\geq\min\{\psi_3(b_j),\psi_3(b_{j-1})\}+O(1/b_j)$ provided $\psi_3(b_{j-1})<\nu_1$.
\end{lemm}

Proof: (i) Suppose $q>q'$ are in the interval. When $P_3$ is decreasing, then  $P_3$ has slope $-1$, hence $P_3(q)-P_3(q')=-(q-q')$, so that
$$\psi_3(q)=\psi_3(q')q'/q-(1-q'/q)=\psi_3(q')-(1+\psi_3(q'))(1-q'/q)<\psi_3(q')$$
by $(3.2)$.\\
When $P_3$ is increasing with slope $\nu_2$, then $P_1(q)-P_1(q')\geq\nu_2(q-q')$, yielding
$$\psi_3(q)\geq\psi_3(q')q'/q+\nu_2(1-q'/q)=\psi_3(q')+(\nu_2-\psi_3(q'))(1-q'/q)>\psi_3(q')$$ 
in view of $(3.2)$.\\ 
When $P_3$ is increasing with slope $\nu_1$, then the same computation as above with $\nu_1$ in place of $\nu_2$ yields the desired result since in this case $\nu_1-\psi_3(q')\geq0$ as a consequence of the additional assumption $(3.7)$.

(ii) If $T_j$ is of type 1, then we have
\begin{eqnarray*}
P_3(q)&\geq&\frac{1}{2}\left(P_3(q)+P_2(q)\right)+O(1)\\
&\geq&\frac{1}{2}\left(P_3(b_j)+P_2(b_j)\right)+\frac{1-\nu_2}{2}(b_j-q)+O(1)\\
&=&P_3(b_j)+\frac{1-\nu_2}{2}(b_j-q)+O(1)\\
&=&P_3(b_j)+\frac{\nu_1}{2}(b_j-q)+O(1)
\end{eqnarray*}
hence the assumption $\psi_3(b_j)\geq0$ yields
\begin{eqnarray*}
\psi_3(q)&\geq&\psi_3(b_j)b_j/q+\frac{\nu_1}{2}(b_j/q-1)+O(1/q)\\
&=&\psi_3(b_j)+(\frac{\nu_1}{2}+\psi_3(b_j))(b_j/q-1)+O(1/q)\\
&\geq&\psi_3(b_j)+O(1/b_j).
\end{eqnarray*}

If $T_j$ is of type 2 then 
\begin{eqnarray*}
P_3(q)&\geq&\frac{1}{2}\left(P_3(q)+P_2(q)\right)+O(1)\\
&\geq&\min\left\{\!\frac{1}{2}\left(P_3(a_j)\!+\!P_2(a_j)\right)+\frac{1}{2}(q\!-\!a_j),\frac{1}{2}\left(P_3(b_j)\!+\!P_2(b_j)\right)+\!\frac{1\!-\!\nu_2}{2}(b_j\!-\!q)\!\right\}+O(1)\\
&=&\min\left\{P_3(a_j)+\frac{1}{2}(q-a_j),P_3(b_j)+\frac{1-\nu_2}{2}(b_j-q)\right\}+O(1)\\
&=&\min\left\{P_3(a_j)+\frac{1}{2}(q-a_j),P_3(b_j)+\frac{\nu_1}{2}(b_j-q)\right\}+O(1).
\end{eqnarray*}
If this minimum is $P_3(b_j)+\frac{\nu_1}{2}(b_j-q)$, then we conclude that $\psi_3(q)\geq\psi_3(b_{j})+O(1/b_{j})$ as for type 1 intervals. If the minimum is $P_3(a_j)+\frac{1}{2}(q-a_j)$, the fact that $P_3$ is increasing in $[b_{j-1},a_{j}]$ yields $P_3(a_j)\geq P_3(b_{j-1})+\nu_1(a_j-b_{j-1})$ and hence $\nu_1\leq1/2$ implies
$$P_3(q)\geq P_3(b_{j-1})+\nu_1(q-b_{j-1})+O(1).\leqno(3.8)$$
As $\psi_3(b_{j-1})<\nu_1$ we may thus write
\begin{eqnarray*}
\psi_3(q)&\geq&\psi_3(b_{j-1})b_{j-1}/q+\nu_1(1-b_{j-1}/q)+O(1/q)\\
&=&\psi_3(b_{j-1})+(\nu_1-\psi_3(b_{j-1}))(1-b_{j-1}/q)+O(1/q)\\
&\geq&\psi_3(b_{j-1})+O(1/b_{j-1}).
\end{eqnarray*}

\begin{prop}
Assume that $\op_3^{\bw}<\nu_1$. Then $\up_3^{\bw}=\liminf_{j\rightarrow\infty}\psi_3(b_j)$.
\end{prop}

Proof: Note that $\op_3^{\bw}<\nu_1$ implies $(3.7)$ for $q\geq q_0$ and that $[q_0,\infty)\subseteq\bigcup_{j>j_0}[b_j,b_{j+1}]$, where $b_{j_0-1}\leq q_0\leq b_{j_0}$. To prove the claim, it suffices to show that 
$$\min_{q\in[b_j,b_{j+1}]}\psi_3(q)\in\{\psi_3(b_j),\psi_3(b_{j+1})\}\mbox{ for }j\geq j_0.\leqno(3.9)$$
We will consider the decomposition of $[b_j,b_{j+1}]$ in the union $[b_j,a_{j+1}]\cup T_{j+1}$. Then either $P_3$ has a single maximum $p_j$ in $[b_j,a_{j+1}]$ ( if $T_{j+1}$ is of type 1 as in Figure 1) resp. $P_3$ is increasing in this whole interval ( if $T_{j+1}$ is of type 2 as in Figure 3) and Lemma \ref{LI} (i) implies that the minimum of $\psi_3(q)$ for $q\in[b_j,a_{j+1}]$ is attained at $q=b_j$ or at $q=a_{j+1}$ resp. at $q=b_j$. In the type 1 case, \ref{LI} (ii) implies that the minimum of $\psi_3(q)$ for $q\in T_{j+1}$ is attained at $q=b_{j+1}$ whereas in the type 2 case, $(3.7)$ implies in particular that $\psi_3(b_{j-1})<\nu_1$ and thus \ref{LI} (ii) implies that the minimum of $\psi_3(q)$ for $q\in T_{j+1}$ is attained at $q=b_{j+1}$ or $q=a_{j+1}$. In both cases, as $a_{j+1}$ lies in both intervals (!), the minimum cannot be attained at $a_{j+1}$ and $(3.9)$ is established.
\vspace{2mm}

Quite analoguously, we prove the dual result: 

\begin{prop}
Assume that $\op_3^{\bw}<\nu_1$. Then $\op_1^{\bw}=\limsup_{j\rightarrow\infty}\psi_1(a^*_j)$.
\end{prop}

Proof: Note that $\op_3^{\bw}<\nu_1$ implies $\psi_3(q)<\nu_1$ for $q\geq q_0$, hence $\psi_1(q)+\psi_2(q)>-\nu_1$ for such $q$, and thus in particular $P_2(a^*_j)-P_1(a^*_j)=4c$ implies
$$\psi_1(a^*_j)>-\nu_1/2\mbox{ for }a^*_j>q_0.\leqno(3.10)$$ 
Moreover we have $[q_0,\infty)\subseteq\bigcup_{j>j_0}[a^*_j,a^*_{j+1}]$, where $a^*_{j_0-1}\leq q_0\leq a^*_{j_0}$. To prove the claim, it suffices to show that 
$$\max_{q\in[a^*_j,a^*_{j+1}]}\psi_1(q)\in\{\psi_1(a^*_j),\psi_1(a^*_{j+1})\}\mbox{ for }j\geq j_0.\leqno(3.11)$$
We will consider the decomposition of $[a^*_j,a^*_{j+1}]$ in the union $B_j\cup [b^*_j,a^*_{j+1}]$. Then either $P_1$ has a single minimum $p^*_{j+1}$ in $[b^*_j,a^*_{j+1}]$ ( if $B_{j+1}$ is of type 1) resp. $P_1$ is increasing in this whole interval ( if $B_{j+1}$ is of type 2). Lemma \ref{LII} (i) implies that the maximum of $\psi_1(q)$ for $q\in[b^*_j,a^*_{j+1}]$ is attained at $q=b^*_j$ or at $q=a^*_{j+1}$ resp. at $q=a^*_{j+1}$. In view of $(3.10)$ Lemma \ref{LII} (ii) is applicable: in the type 1 case, it implies that the maximum of $\psi_1(q)$ for $q\in B_{j}$ is attained at $q=a^*_{j}$ whereas in the type 2 case, it implies that the maximum of $\psi_1(q)$ for $q\in B_{j}$ is attained at $q=b^*_{j}$ or $q=a^*_{j}$. Here $b^*_{j}$ lies in both intervals, so in both cases the maximum cannot be attained at $b^*_{j}$ and $(3.11)$ is established.
\vspace{2mm}

\section{A Jarnik type relation between $\op_3^{\bw}$ and $\up_1^{\bw}$}

\begin{theo}
Assume that $\op_3^{\bw}<\nu_1$. Then 
$$(2-\nu_1)\op_1^{\bw}+(1+\nu_1)\up_3^{\bw}\leq-3\op_1^{\bw}\up_3^{\bw}\leq(2-\nu_2)\op_1^{\bw}+(1+\nu_2)\up_3^{\bw}.\leqno(4.1)$$
\end{theo}

Proof: From $(3.6)$ in Lemma \ref{L+L} we have
$$\psi_1(a^*_j)=-\frac{(1+\nu(j))\psi_3(b_j)}{3\psi_3(b_j)+(2-\nu(j))}+O(1/b_j)=:-f_{\nu(j)}(\psi_3(b_j))+O(1/b_j),\leqno(4.2)$$
where $\nu_1\leq\nu(j)\leq\nu_2$ depending on the position of $r_j$ in the interval $I_j$. By Proposition 3.5 and the fact that $1/b_j\rightarrow0$ for $j\rightarrow\infty$ we obtain
\begin{eqnarray*}
\op_1^{\bw}=\limsup_{j\rightarrow\infty}\psi_1(a^*_j)&=&\limsup_{j\rightarrow\infty}-f_{\nu(j)}(\psi_3(b_j))\\
                      &=&-\liminf_{j\rightarrow\infty} f_{\nu(j)}(\psi_3(b_j)).
\end{eqnarray*}

Note that $f_{\nu(j)}(x)$ is increasing for $x>(\nu-2)/3$, hence for all possible values of $\psi_3(b_j)$ and that for fixed $x$ in this range $\nu<\mu$ implies $f_{\nu}(x)<f_{\mu}(x)$. We may thus write
\begin{eqnarray*}
\liminf_{j\rightarrow\infty} f_{\nu(j)}(\psi_3(b_j)&\geq&\min_{j}f_{\nu(j)}(\liminf_{j\rightarrow\infty}\psi_3(b_j))\\
&=&\min_{j}f_{\nu(j)}(\up_3^{\nu})\\
&\geq&f_{\nu_1}(\up_3^{\bw}).
\end{eqnarray*}
Alltogether we have 
$$\op_1^{\bw}\leq-f_{\nu_1}(\up_3^{\bw}),$$
which gives exactly the left hand inequality of $(4.1)$.

On the other hand $(3.6)$ may be written as
$$\psi_3(b_j)=-\frac{(2-\nu(j))\psi_1(a^*_j)}{3\psi_1(a^*_j)+(1+\nu(j))}+O(1/b_j)=:-g_{\nu(j)}(\psi_1(a^*_j))+O(1/b_j)\leqno(4.3)$$
and by Proposition 3.4

\begin{eqnarray*}
\up_3^{\bw}=\liminf_{j\rightarrow\infty}\psi_3(b_j)&=&\liminf_{j\rightarrow\infty}-g_{\nu(j)}(\psi_1(a^*_j))\\
                      &=&-\limsup_{j\rightarrow\infty} g_{\nu(j)}(\psi_1(a^*_j)).
\end{eqnarray*}
Like before, $g_{\nu}(x)$ is increasing for $x>-\frac{1+\nu}{3}$ and as in Proposition 3.5 we use the fact that $\op_3^{\bw}<\nu_1$ implies $\psi_1(a^*_j)>-\nu_1/2\geq-1/4>-1/3$ for $j\geq j_0$, to conclude that $g_{\nu}(x)$ is increasing for all values of $\psi_1(a^*_j)$ where $j\geq j_0$.
Still for fixed $x$ in this range $\nu<\mu$ implies $g_{\nu}(x)<g_{\mu}(x)$. This yields
\begin{eqnarray*}
\limsup_{j\rightarrow\infty} g_{\nu(j)}(\psi_1(a^*_j)&\leq&\max_{j}g_{\nu(j)}(\limsup_{j\rightarrow\infty}\psi_1(a^*_j))\\
&=&\max_{j}g_{\nu(j)}(\op_1^{\bw})\\
&\leq&g_{\nu_2}(\op_1^{\bw}).
\end{eqnarray*}
Alltogether we have 
$$\up_3^{\bw}\geq-g_{\nu_2}(\op_1^{\bw})$$
which gives exactly the right hand inequality of $(4.1)$.
\vspace{2mm}

{\bf{Remark}}: For sake of symmetry, we may write $(4.1)$ as
$$(1+\nu_2)\op_1^{\bw}+(1+\nu_1)\up_3^{\bw}\leq-3\op_1^{\bw}\up_3^{\bw}\leq(1+\nu_1)\op_1^{\bw}+(1+\nu_2)\up_3^{\bw}.$$
Note that for $\nu_1=\nu_2=1/2$ both bounds for $-3\op_1^{\bw}\up_3^{\bw}$ are identical and the quantities $\op_1^{\bw},\up_3^{\bw}$ are those from the classical setting. Hence writing them as $\op_1,\up_3$ and multiplying by $2/3$ yields the classical Jarnik identity
$$\op_1+\up_3+2\op_1\up_3=0.$$
The assumption $\op_3<1/2$ can be dropped in this case.

\vspace{15mm}

\begin{center}{\Large{\bf References}}\end{center}

\vspace{4mm}

[1] S. Chow, A. Ghosh, L. Guan, A. Marnat, D. Simmons; {\emph{Diophantine transference inequalities: weighted, inhomogeneous, and intermediate exponents}}. Annali Della Scuola Normale Superiore Di Pisa, arXiv:1808.07184v2 
\vspace{4mm}

[2] O. German; {\emph {Transference theorems for Diophantine approximation with weights}}. https://arxiv.org/abs/1905.01512v2
\vspace{4mm}

[3]  A. Marnat; {\emph { About Jarnik's-type relation in higher dimension}}. Ann. Inst. Fourier (Grenoble) {\bf 68} No. 1 (2018), p. 131-150 
\vspace{4mm}

[4] D. Roy; {\emph{ On the topology of Diophantine approximation spectra}}. Compos. Math. {\bf 153} No. 7 (2017), p. 1512-1546
\vspace{4mm}

[5] W.M. Schmidt; {\emph{On Parametric Geometry of Numbers}}. Submitted at Acta Arithmetica.
\vspace{4mm}

[6] W. M. Schmidt, L. Summerer; {\emph{Parametric Geometry of Numbers and applications}}. Acta Arithmetica  {\bf 140} No. 1 (2009), p. 67-91 
\vspace{4mm}

[7] W. M. Schmidt, L. Summerer; {\emph{The generalization of Jarnik's identity}}. Acta Arith. {\bf 175} No. 2 (2016), p. 119-136
\vspace{4mm}

[8] L. Summerer; {\emph{Generalized simultaneous approximation to $m$ linearly dependent reals}}. Mosc. J. Comb. Number Theory {\bf 8} No. 3 (2019), p. 219-228

\vspace{8mm}

\noindent
Leonhard Summerer\\
Fakult\"at f\"ur Mathematik\\
Universit\"at Wien\\
Oskar-Morgenstern-Platz 1\\
A-1090 Wien, AUSTRIA\\
E-mail: leonhard.summerer@univie.ac.at

\end{document}